\documentclass{article}

\usepackage{amsmath}
\usepackage{amsthm}
\usepackage{amsfonts}
\usepackage{amssymb}
\usepackage{amscd}
\usepackage{diagrams}

\newtheorem{thm}{Theorem}
\newtheorem{defi}[thm]{Definition}
\newtheorem{prop}[thm]{Proposition}
\newtheorem{cor}[thm]{Corollary}
\newtheorem{lem}[thm]{Lemma}
\newtheorem{ex}[thm]{Example}

\newcommand{\io}{\iota}

\newcommand{\Z}{\mathbb Z}
\newcommand{\Q}{\mathbb Q}
\newcommand{\R}{\mathbb R}
\newcommand{\N}{\mathbb N}

\newcommand{\cp}{\mathbb{CP}}

\newcommand{\Di}{\mathcal D}

\newcommand{\ot}{\otimes}

\newcommand{\de}{\partial}

\def\z2{\Z_2}

\newarrow{Equal} =====
\newarrow{Hook} C--->

\def\Om{\Omega}
\def\Si{\Sigma}
\def\x{\times}
\def\sm{\wedge}
\def\La{\Lambda}
\def \osn{\Om^n \Si^n X}

\begin{document}

\author{Paolo Salvatore}
\title{Configuration spaces on the sphere and higher loop spaces}

\date{}

\maketitle

\begin{abstract}

We show that the homology over a field of the
space $\La^n \Si^n X$ of free maps from the $n$-sphere to
the $n$-fold suspension of $X$ depends only on
the cohomology algebra of $X$ and compute it explicitly.
We compute also the homology of the closely related labelled
configuration space
$C(S^n,X)$ on the $n$-sphere with labels in $X$
 and of its completion,
that depends only on the homology of $X$.
In many but not all cases the homology
of $C(S^n,X)$
coincides with the homology of
$\La^n \Si^n X$.
In particular we obtain the homology
of the unordered configuration spaces on a sphere.

\


\end{abstract}

\section*{Introduction}

In order to compute the homology of
the space of based maps $\Om^n(Y)=map_*(S^n,Y)$ from the
$n$-sphere to a CW-complex $Y$, in general one needs to know a great
deal of information on $Y$.
However, when $Y=\Si^n X$ is a $n$-fold suspension,
a classical result by Milgram \cite{Milgram}
 states that the homology of
$\Om^n \Si^n X$ depends just on the homology
of $X$. An explicit description is given in \cite{Cohen}.
This depends upon the existence of a small model
for the mapping space.
Namely the configuration space $C(\R^n,X)$ of pairwise
distinct points in $\R^n$ with labels in $X$, modulo base point
cancellation,
is homotopy equivalent to
$\osn$ when $X$ is connected \cite{May}.
In general $\osn$ is the group completion of $C(\R^n,X)$.

Let us turn our attention to the free mapping space $\La^n Y =
map(S^n,Y)$. In the case $n=1$ the homology of $\La \Si X$ depends
only on the homology of $X$ \cite{Menichi}. This depends again on
the existence of a small model for the mapping space, the
configuration space $C(S^1,X)$ of points in the circle with labels
in $X$ \cite{Bodig}, for $X$ connected. This does not extend to
all $n$. For $X$ connected the configuration space $C(S^n,X)$ is
homotopy equivalent to the section space of a bundle over $S^n$
with fiber $\Si^n X$ \cite{Bodig}, obtained from the tangent
bundle of the sphere $\tau_n$ by adding a point at infinity to
each fiber and smashing it with $X$. In general this section space
is a kind of completion of $C(S^n,X)$ . When $n \neq 1,3,7, \,
S^n$ is not parallelizable, so that one does not expect $C(S^n,X)$
to be homotopy equivalent to the free mapping space
$\Lambda^n\Si^n X$.

\smallskip

In this paper we compare
the homotopy types of these spaces
and we compute their homology with coefficients in any field.
Up to homotopy there are fibrations
$C(\R^n,X) \to C(S^n,X) \to \Si^n X$ (for $X$ connected)
and
$\osn \to \La^n \Si^n X \to \Si^n X$,
induced by the evaluation at one point.
Since the basis of both fibrations
is a suspension, we can
reconstruct
the total spaces
by means of the clutching functions $\Si^{n-1}X \x \osn \to \osn$.
It turns out that the clutching
function of the first fibration
 comes from
the action of the little $n$-discs operad
$S^{n-1} \x X \x C(\R^n,X) \to C(\R^n,X)$.
This action adds to a
`cloud' of labelled points in
$C(\R^n,X)$ an extra point
in the direction parametrized by $S^{n-1}$ with label
parametrized by $X$.
The same argument works when $X$ is not connected.
In general we get decompositions both for
$C(S^n,X)$ (Theorem \ref{conf}) and its completion
(Theorem \ref{sec}),
that coincide for $X$ connected.
The decomposition is compatible with the Snaith
splitting (Proposition \ref{stable}).
The clutching function induces in homology the so called
Browder operation, that has been computed explicitly
for $C(\R^n,X)$ in \cite{Cohen}.
For example in characteristic 0 and for $n=2$ we get the
adjoint action of the homology of $X$
on the free Gerstenhaber algebra
that it generates.
This allows to determine the homology
of $C(S^n,X)$, that depends only on the homology of $X$
(Corollary \ref{corconf}).
In particular for $X=S^0$ we obtain the homology
of the unordered configuration spaces on $S^n$.
We write down an explicit basis of their homology in
Theorem \ref{basisconf}.
The case when $n$ is even and the characteristic is odd
is not covered by the methods
in \cite{BCT}.

Let us consider now the evaluation fibration of the mapping space
$\La^n \Si^n X$. In this case the clutching function is obtained
by twisting the first argument $S^{n-1} \x X \to \osn$  and  then
applying as before the little $n$-discs action $S^{n-1} \x
(\osn)^2 \to \osn$ (Theorem \ref{maps}) . The twist is the
composition operation induced by the element $c \in
\pi_{n-1}(\Om^n_1 S^n)$ adjoint to the Whitehead product of the
generator $\pi_n(S^n)$ with itself. This happens exactly because
$c$ classifies the fiberwise compactification of the tangent
bundle of the $n$-sphere (Lemma \ref{lemmino}).

Now $c$ goes in homology
(up to shift of component)
to the Browder operation
of the identity $\io \in H_0(\Om^n S^n)$ with itself.
This implies that the clutching function of the mapping
space in homology
depends only on the
cohomology algebra of $X$ (for finite type),
and so does the homology of $\La^n \Si^n X$ (Corollary \ref{cormaps}).
We show in
 Example \ref{mod3}
that the homology of $\La^n \Si^n X$ does not depend just on the
additive homology of $X$. As an application we write down an
explicit basis of the homology of the mapping space $\La^n S^m$,
for $n \leq m$, in Theorems \ref{basisconf} and \ref{basismaps}.

\smallskip

Let us compare the configuration space
(or its completion for $X$ not connected)
 and the mapping space:
when $c$ vanishes they are homotopy equivalent. This happens when
either $X$ is a suspension, or $n=1,3,7$, or for $n$ odd after
inverting the prime 2 (Proposition \ref{equal}). Moreover the two
spaces have the same homology either mod 2 (Corollary \ref{mod2}),
or for $X$ connected and rationally (Corollary \ref{mod0}). In
general they do not have the same homology for $n$ even at odd
primes (Example \ref{mod3}), or rationally when $X$ is not
connected (Example \ref{not0}).

We present also a homotopy pullback decomposition of our
spaces, induced by decomposing the $n$-sphere of the domain
as union of two discs (Propositions \ref{homaps} and \ref{hosec}).
In this case the twist occurs in the diagram
describing the configuration space.
We apply the induced sequence of homotopy groups
to show in Example
\ref{last}
 that the configuration space and the mapping space
are not rationally homotopic
for $n=2$ and $X=\cp^2$, but they
have the same homology over
all fields.

It would be interesting to phrase our computations
in terms of the higher Hochschild homology spectral
sequence described in 5.9 of
 \cite{Pirashvili}. 

\

In sections \ref{2},\ref{3} and \ref{4}
 we consider respectively
labelled configuration spaces on the sphere,
the associated spaces of sections, and mapping spaces.
In section \ref{5} we compute some examples and
we compare the homology of the spaces above.
In section \ref{6} we present the homotopy pullback decomposition
and an application.

\

I am grateful to Sadok Kallel for many helpful conversations.

\bigskip

\section{Configuration spaces} \label{2}

Let $X$ be a based CW-complex, not necessarily connected.
Let $n$ be a positive integer.

We recall some definitions from \cite{Bodig}.
Let $F_k(M)$ be the space of
pairwise distinct $k$-tuples in a manifold $M$.
We consider the elements of the space  $F_k(M) \x_{\Si_k} X^k$
as finite sets of points in $M$ with a label in $X$.
The configuration space $C(M,X)$ is the quotient of
the disjoint union
$\coprod_k  F_k(M) \x_{\Si_k} X^k$ under the equivalence relation
induced by cancelling points
labelled by the base point of $X$.
Given a closed subset $A \subset M$, the relative configuration space
$C(M,A;X)$ is the quotient induced by cancelling
all points located in $A$ and those labelled by the base point.

We identify $S^n =\R^n \cup \infty$, based at infinity.
Let $D_n \subset \R^n$ be the unit disc.
The sequence $D_n \to S^n \to (S^n,D_n)$,
for $X$ connected,
induces a quasifibration
$C(D_n,X) \to C(S^n,X) \to C(S^n,D_n;X)$
(Lemma p.178 in \cite{Bodig}).
Now $\Si^n X \subset C(S^n,D_n;X)$ is a strong deformation retract,
seen as the subspace of configurations with a single labelled point
\cite{McDuff}.
We will construct a homotopy pushout decomposition of $C(S^n,X)$,
by pulling back to the total space the decomposition of
$\Si^n X = S^n \sm X$ induced by splitting
$S^n$
as union of the $n$-discs $D_n$ and $E_n= \overline{S^n -D_n}$.
The quasifibration over $D_{n+} \sm X$ and $E_n \sm X$ is homotopically
trivial. It remains to glue the two pieces together.
This idea works even when $X$ is not connected.

We recall that $C(D_n,X) \simeq C(\R^n,X)$ is homotopy equivalent to
the free algebra
on the operad of little $n$-discs generated by $X$.
In particular there is an action of the space of two little $n$-discs
$\Di_n(2)$ by a map
$\rho_2: \Di_n(2) \x C(\R^n,X) \x C(\R^n,X) \to C(\R^n,X)$.
Let us restrict this map to the deformation retract
 $S^{n-1} \subset \Di_n(2)$ in the first factor
and to  $X \subset C(\R^n,X)$, the subspace
of configurations consisting of a single labelled point
located at the origin, in the second
factor.
We get a map
$\lambda: S^{n-1}_+ \sm X \x C(\R^n,X) \to C(\R^n,X)$, that we shall
call the {\em Browder map}.
Up to deformation
$\lambda(z,x,c)$ is the configuration
obtained from the configuration $c$ by adding
a point outside
$c$ in direction $z$ with label $x$.

\begin{thm} \label{conf}
There is a homotopy pushout
$$\begin{diagram}
S^{n-1}_+ \sm X \x C(\R^n,X) & \rTo^\pi & X \x  C(\R^n,X) \\
\dTo^\lambda                     &        &  \dTo        \\
C(\R^n,X)                & \rTo            &  C(S^n,X) ,
\end{diagram}$$
where $\lambda$ is the Browder map and
$\pi$ the projection.
\end{thm}

\begin{proof}

Let us consider the closed deformation retract $K$ of
$C(S^n,X)$ containing the configurations
with at most one labelled point located outside $D_n$.
The deformation is obtained by pushing the particles
away from the north pole and deforming slightly the labels
towards the base point.
Let
$A,B \subset K$ be the closed subsets
such that the point located outside $D_n$, if it exists,
is respectively in $2D_n$ and $2E_n$.
Clearly $A \cup B=K$.
Moreover the inclusion $C(D_n,X) \subset A$ is a homotopy equivalence.
On the other hand
$A \cap B \cong   S^{n-1}_+\sm X \x C(D_n,X)$,
and $B \cong {D_n}_+ \sm X \x C(D_n,X)$.
With these identifications the map $A \cap B \to B$ is induced by the inclusion $S^{n-1} \subset D_n$.
Moreover the inclusion
$A \cap B \to A$ is homotopic to
the Browder map.

\end{proof}

From now on we consider homology with coefficients in a field $K$.
We recall the definition of the Browder operation:
let $Y$ be a space acted on by the little
$n$-discs operad.
Let $e_{n-1} \in H_{n-1}(\Di_n(2))$ be the class represented by the
deformation retract $S^{n-1} \subset \Di_n(2)$.
Then the Browder operation,
           for  $x \in H_p(Y)$ and $y \in H_q(Y)$,
is defined by $[x,y]:= \rho_{2*}(e_{n-1} \ot x \ot y)
\in H_{p+q+n-1}(Y)$. This differs from (5.7) of \cite{Cohen}
by the sign $(-1)^{(n-1)p+1}$.

For $n=1$ the fundamental class of $S^0$ is $\{1\}-\{-1\}$, so that
the Browder operation is the commutator of the Pontrjagin product.

\begin{defi}
The Browder action is the homomorphism

$$[\;,\;]=\lambda_*:
\Si^{n-1}\tilde{H}(X) \ot H(C(\R^n,X)) \to H(C(\R^n,X)).$$
\end{defi}

Recall from Theorem 3.1 in \cite{Cohen}
that the homology
$\tilde{H}( C(\R^n,X))$ depends functorially on the homology $\tilde{H}(X)$.
For example in characteristic 0 and for $n>1$
the former is the  free $n$-algebra
on the latter.
A $n$-algebra is
a graded variation of a Poisson algebra \cite{nat}. In our case
the Pontrjagin product
 is the commutative product and the Browder operation is the
bracket (up to sign). The Browder action for $char(K)=0$ is
exactly the adjoint action of the space of  generators
$\tilde{H}(X)$ on  $H(C(\R^n,X))$.

In positive characteristic
there are additional homology operations generating
$H(C(\R^n,X))$ over $H(X)$,
the unstable Dyer-Lashof operations.
See (1.1),(1.3) of \cite{Cohen} and Proposition \ref{dyla}.

In general the Browder action can be computed explicitly
and depends just on the homology of $X$
by Thm. 1.2 (5,6,8) and Thm. 1.3 (4) of \cite{Cohen}.
For example mod 2 the action is trivial
 on all Dyer-Lashof operations
except the top one.
Recall that $H(C(\R,X))$ is the free associative algebra on
$\tilde{H}(X)$.

The Mayer-Vietoris sequence implies the following corollary.

\begin{cor}  \label{corconf}
The homology of $C(S^n,X)$ depends only on the homology of $X$ and
is the direct sum of the suspended kernel and the cokernel of the
Browder action on the homology of $C(\R^n,X)$.
\end{cor}

The configuration spaces $C(M,X)$, for a manifold $M$,
are naturally filtered by subspaces $C_k(M,X)$, having
at most $k$ particles.
We denote the subquotients by $D_k(M,X)=C_k(M,X)/C_{k-1}(M,X)$
and $D_k^n(X)=D_k(\R^n,X)$.
For example the space $D_k(S^n,S^0)$ is the
unordered
configuration space on the sphere
$C_k(S^n)=F_k(S^n)/\Si_k$.

The homology
$H(C(M,X))$ splits naturally as
sum of $H(D_k(M,X))$, as $k$ varies.
Compare Lemma 4.2 p.238 in \cite{Cohen}.
This explains the first assertion of Corollary \ref{corconf},
as

$H(C(S^n,X)) \cong \bigoplus_k H(F_k(S^n); \tilde{H}(X)^{\ot k} ) .$

The Browder map is compatible with the
filtration by number of particles,
so that we have the following refinement:

\begin{cor} \label{split}
The reduced homology of $D_k(S^n,X)$, for $k>1$,
is the direct sum of the suspended kernel and the cokernel
of the Browder action
$\Si^{n-1}\tilde{H}(X) \ot \tilde{H}(D_{k-1}^n(X)) \to
\tilde{H}(D_k^n(X))$.
\end{cor}

There is a geometric version of this corollary.
Recall that $C(M,X)$ is stably homotopy equivalent
to the wedge sum of
$D_k(M,X)$, indexed over $k$ \cite{Bodig}.
The Browder map defines a stable map
$b_k: \Si^{n-1}X \sm D_{k-1}(X) \to D_k(X)$,
by restriction to the relevant wedge summands.
For example
the map $b_2$ is the restriction to the stable
summand $S^{n-1} \sm X \sm X$
of the quotient map
$S^{n-1}_+ \sm X \sm X \to S^{n-1}_+\sm_{\Z_2}(X \sm X)$.

\begin{prop}              \label{stable}
The $k$-th summand of the stable splitting
$D_k(S^n,X)$ is stably homotopy equivalent to
the cofiber of
$b_k: \Si^{n-1}X \sm D_{k-1}^n(X) \to D_k^n(X)$ for $k >1$ and to
$\Si^n X \vee X$ for $k=1$.
\end{prop}

\begin{proof}

We obtain a homotopy pushout $(k)$ expressing the homotopy type of
$C_k(S^n,X)$, if we replace, in the diagram of Theorem 1, $C$ by
$C_{k-1}$ in the top row and by $C_k$ in the bottom row. For $k=1$
this proves the proposition. For $k >1$ the cofiber diagram of
$(k-1)$ and $(k)$ gives $D_k(S^n,X)$ as homotopy pushout. The top
map $\pi_k$ of such diagram is a stable retraction, and $b_k$ is
exactly the restriction of the left hand side map $\lambda_k$ to
the stable kernel of $\pi_k$.
\end{proof}

\medskip

\section{Section spaces} \label{3}

We recall that there is a `scanning' map
$C(D_n,X) \to \osn \cong map_*(D_n/S^{n-1},D_n/S^{n-1} \sm X)$.
Let $i:X \to \osn$ be the adjoint of the identity.
Roughly speaking, if we scan a configuration,
we obtain a $n$-fold loop, sending small
discs centered at the configuration points to
$\Si^n X \cong D_n/S^{n-1} \sm X$, by
identifying them to the unit disc and applying
$i$ to the label of the center.
The $n$-fold loop sends everything else to the base point.

The scanning map can be extended to $C(S^n,X)$ by using the
exponential map of the sphere. Let $\Si^n X \to \tau_n^+ ( X )\to
S^n$ be the bundle obtained by smashing $X$ fiberwise with the
fiberwise one-point compactification $\tau_n^+$ of the tangent
bundle $\tau_n$ of the $n$-sphere. Then the scanning map goes from
$C(S^n,X)$ to the space $\Gamma(\tau_n^+ ( X))$
 of sections
of the bundle.
Up to homotopy the scanning maps fit into a diagram
from the sequence

$$C(D_n,X) \to C(S^n,X) \to C(S^n,D_n;X)$$ to the
evaluation fibration
$$\osn \to \Gamma(\tau_n^+ ( X)) \rTo^{ev} \Si^n X.$$

The scanning maps $C(D_n,X) \to \osn$ and $C(S^n,X) \to
\Gamma(\tau_n^+ ( X))$ are homotopy equivalences if and only if
$X$ is connected.

We can define the Browder map
 for the $n$-fold loop
space $\osn$, similarly as for the configuration space, by
using the action
$\Di_n(2) \x \osn \x \osn \to \osn$
 of the $n$-discs operad,
and restricting the action in the second factor to
the subspace $X \rTo^i \osn$.

\

If $X$ is connected then
the next theorem reduces to Theorem \ref{conf}.
Also the proof is an adaptation of the proof of Theorem \ref{conf}
to the general case.

\begin{thm} \label{sec}
There is a homotopy pushout
$$\begin{diagram}
S^{n-1}_+ \sm X \x \osn & \rTo^\pi & X \x  \osn \\
\dTo^\lambda                     &        &  \dTo        \\
\osn                & \rTo            & \Gamma(\tau_n^+ ( X))   ,
\end{diagram}$$
where $\lambda$ is the Browder map and
$\pi$ the projection.
\end{thm}

\begin{proof}
Let $ev: \Gamma(\tau_n^+ ( X)) \to \tau_n(N)^+ \sm X \cong \Si^n
X$ be the fibration evaluating at the north pole $N$.

We construct a map $(\alpha,\beta,\gamma)$ from the diagram
$$D_{n+} \sm X \x \Om^n\Si^n X \lTo S^{n-1}_+ \sm X \x \Om^n\Si^n X \rTo M(\lambda),$$
where $M(\lambda)$ is the mapping cylinder of the Browder map,
to the diagram
$$ev^{-1}(D_{n+} \sm X) \lTo ev^{-1}(S^{n-1}_+ \sm X) \rTo ev^{-1}(E_n\sm X).$$

We construct $\alpha$ so that it covers an automorphism of $D_{n+}
\sm X$ homotopic to $-1 \sm X$.
 We identify both the southern and the northern
hemisphere to the unit disc, by projection from the opposite pole.
This gives a trivialization of the tangent bundle on each
hemisphere.

The map $\alpha(z \sm x,y) \in \Gamma(\tau_n^+ ( X))$ is defined
by $y \in Map(D_n,\de D_n; \Si^n X)$ on the southern hemisphere.
Let $\phi:(D_n,S^{n-1}) \to (S^n,\infty)$ be the relative
homeomorphism $\phi(x)=x/(1-|x|)$.
 In the northern hemisphere we have the adjoint $n$-fold loop
$\alpha(z \sm x,y)(w)= \phi(2(w-z/4)) \sm x$ if $|w-z/4| \leq
1/2$, and otherwise $w$ goes to the point at infinity.

Similarly, on the cylinder $I \x S^{n-1}_+ \sm X \x \Om^n\Si^n X
\subset M(\lambda)$, $\gamma(t,z \sm x,y)$ is defined in the
southern hemisphere by $y$. On the northern hemisphere $\gamma(t,z
\sm x,y)(w)= \phi(2(w-z(t+1)/4)) \sm x$ for $|w-z(t+1)/4| \leq
1/2$ and otherwise $w$ goes to the point at infinity. The map
$\gamma$ sends the end of the mapping cylinder $\Om^n \Si^n X$ to
$ev^{-1}(*)$ by projection from the north pole.

Both $\alpha$ and $\beta$ are homotopy equivalences, because they
cover an automorphism of the basis and  induce homotopy
equivalences of the fibers over points of all components. Finally
$\gamma$ is a homotopy equivalence, since the domain has $\osn$ as
deformation retract, identified by $\gamma$ to $ev^{-1}(*) \simeq
ev^{-1}(E_n \sm X)$. By the homotopy invariance of homotopy
pushouts the theorem is proved.

\end{proof}

Also in this case we call Browder action
the homomorphism
$\lambda_*:\Si^{n-1} \tilde{H}(X) \ot H(\osn) \to H(\osn)$
induced by the Browder map.

\begin{cor}       \label{corsec}
The homology of $\Gamma(\tau_n^+ ( X))$ depends only on the
homology of $X$, and is the direct sum of the suspended kernel and
the cokernel of the Browder action on the homology of $\osn$.
\end{cor}

\begin{proof}

Recall that $H(\osn)$ is obtained from $H(C(\R^n,X))$,
according to the
group completion theorem,
by adding the inverses of the components.

In the case $n>1$
$\pi_0(C(\R^n,X)) = \N[\pi_0(X)]$, and
$\pi_0(\osn) = \Z[\pi_0(X)]$.
We must only compute the Browder action on the inverse
of a component $a \in \pi_0(X) \subset H(\osn)$.
But by the Poisson relation and since the Browder action
on the unit is trivial,
we obtain $[x,a^{-1}]=- [x,a]a^{-2}$.

In the case $n=1$ $\pi_0(C(\R,X))$ and
$\pi_0(\Omega \Sigma X)$ are respectively the free monoid and
the free group
on the based set $\pi_0(X)$, and the Pontrjagin product is clear.

\end{proof}

\begin{ex} \label{massa}
If $X=S^0 \vee S^0$, then $\Gamma(\tau_n^+  (X))$ is homotopy
equivalent to the section space of the bundle $F_3(\R^{n+1}) \to
F_2(\R^{n+1})$ defined by forgetting the third point of a
configuration {\rm \cite{Massey}}.
\end{ex}

\section{Higher loop spaces} \label{4}

In some special cases the section space and the mapping space
coincide.

\begin{prop}   \label{equal}

The fibration $\Si^n X \to \tau_n^+ ( X) \to S^n$ is trivial, so
that the section space is homotopic to the ordinary mapping space
$\La^n \Si^n X$, in the following cases:
\begin{enumerate}
\item When $X$ is a suspension $X=\Si Y$;

\item For $n=1,3,7$;

\item  For $n$ odd and away from the prime 2.
\end{enumerate}
\end{prop}

\begin{proof}

1) The tangent bundle of the sphere is trivialized by adding a
trivial line bundle. Thus $\tau_n^+ ( \Si Y) \cong (\tau_n \oplus
\R)^+ ( Y ) \cong \Si^{n+1} Y \x S^n$.

2) The spheres $S^1, S^3, S^7$ are parallelizable.

3) Follows from the following lemma, since
$2[\iota_n,\iota_n]=0 \in \pi_{2n-1}(S^n)$ for $n$ odd,
 where $\io_n \in \pi_n(S^n)$ is the generator
and the brackets denote the Whitehead product.
\end{proof}

\begin{lem}      \label{lemmino}
The fibration $S^n \to \tau_n^+ \to S^n$ with section $S^n \to
\tau_n^+$ at infinity is classified by the adjoint of the
Whitehead product $[\iota_n,\iota_n] $ in $\pi_{n-1}(\Omega^n _1
S^n)$.
\end{lem}
\begin{proof}
The bundle $\tau_n^+$ is trivial if we forget the section, because
it is the sphere bundle of $\tau_n \oplus \R$, and its section can
be identified to the diagonal $\Delta:S^n \to S^n \x S^n$.

We will show that the homotopy class of the clutching function
$c:S^{n-1} \to \Omega^n_1 S^n$ of the fibration $\tau_n^+$ is
exactly the boundary of $\iota_n \in \pi_n(S^n)$ in the long exact
sequence associated to the evaluation fibration $\Omega^n_1 S^n
\to \Lambda^n_1 S^n \to S^n$. The index 1 denotes the component of
degree 1 maps. We identify $S^{n-1}$ to the equator of $S^n$. We
take the north pole $N$  and the south pole  $S$ respectively as
base points in the domain and the range.

Let $H:I \x S^{n-1} \to SO(n+1) \to \La^n_1 S^n$ be the
transformation such that $H(t,a)$ rotates by $t\pi$ the plane
generated by $a$ and $N$, and fixes the orthogonal complement. The
clutching function of $\tau_n^+$ is $c=H(1,\_):S^{n-1} \to \Om^n_1
S^n$. Let us identify the unreduced suspension of $S^{n-1}$ to
$S^n$ via $[t,a] \mapsto H(t,a)(N)$. Then $H$ defines an element
in $\pi_n(\La^n_1 S^n,\Om^n_1 S^n)$ projecting to the generator
$\io_n \in \pi_n(S^n)$, and $c=\de \io_n$. But by a theorem  of
Whitehead \cite{Wh} $\de(\iota_n)$ is adjoint to
$[\iota_n,\iota_n]$, up to sign convention.
\end{proof}

Let $w: S^{n-1}_+ \sm X \to  \Om^n\Si^n X$ be the composite of
$c_+ \sm i: S^{n-1}_+ \sm X \to \Om^n_1 S^n_+ \sm \Om^n \Si^n X$
and of the composition operation $o:\Om^n_1 S^n_+ \sm \Om^n \Si^n
X  \to \osn$ given by $o(\psi,f)=f \psi$.

We call the map
$k(s, x,y)=\lambda(s,w(s, x),y)$ the {\em twisted} Browder map.

\begin{thm} \label{maps}

There is a homotopy pushout

$$\begin{diagram}
S^{n-1}_+ \sm X \x \Om^n \Si^n X & \rTo^\pi & X \x \Omega^n \Si^n X  \\
\dTo^k                           &        &  \dTo        \\
\Omega^n \Si^n X    & \rTo            &  \La^n\Si^n X ,
\end{diagram}$$

where $k$ is the twisted Browder map
and $\pi$ is the projection.

\end{thm}

\begin{proof}

We identify the two hemispheres
to the unit disc by orientation preserving
diffeomorphisms and we proceed as in Theorem \ref{sec}.
The main difference is that the trivializations
regard the domain but not the range of the mapping space.

Thus the twisted Browder map $k$ replaces $\lambda$ in the proof,
by comparing
the trivializations of the one-point compactified tangent bundle
on the two hemispheres with
Lemma \ref{lemmino}.

\end{proof}

\medskip

In the next corollary $X$ is supposed to have finite type.
Otherwise we must replace the expression `cohomology algebra' by
`homology coalgebra'.

\begin{cor}      \label{cormaps}
The homology of $\La^n \Si^n X$ with coefficients in a field
depends only on the cohomology algebra of $X$
and is the direct sum of the suspended kernel and the cokernel
of the twisted Browder action .
\end{cor}

\begin{proof}
Let $\iota \in H_0(\Om^n_1 S^n)$ be
the class of the identity map.
By the diagram at p. 215 in \cite{Cohen}
the Hurewicz homomorphism sends the
adjoint Whitehead product in $\pi_{n-1}(\Om^n_2 S^n)$ to
the Browder operation $[\iota,\iota] \in H_{n-1}(\Om^n_2 S^n)$,
with our sign convention.
Therefore, for $x \in H(X)$,
$w_*(e_{n-1} \ot x) =
o_*([\io,\io] \io^{-1} \ot x)$ is the composition product
in homology,
where we identify $x$ to $i_*(x) \in H(\osn)$.
Suppose that
the iterated coproduct $\Delta:H(X) \to H(X^3)$ is
$\Delta(x)=\sum_i x_i' \ot x_i'' \ot x_i'''$.

Let $\chi$ be the conjugation of the Hopf algebra
$H(\osn)$. It depends only on the cohomology algebra of $X$,
by induction and by Proposition 1.5 of \cite{Cohen}.

Then the twisted Browder action is
$k_*(x \ot y)= [x,y] + \sum_i[x_i',x_i'']\chi(x_i''')y.$

Compare also
Theorem 3.2 (iv) of \cite{533II}.

\end{proof}

\section{Examples and applications} \label{5}

\begin{cor} \label{mod2}
At the prime 2 the homology of $\La^n\Si^n X $ depends just on the
homology of $X$, and is isomorphic to the homology of the section
space $\Gamma(\tau_n^+ ( X))$.
\end{cor}
\begin{proof}
 The Browder operation
$[\io,\io]$ is trivial mod 2.
\end{proof}

\begin{cor} \label{mod0}
If $X$ is connected then  $\Lambda^n \Sigma^n X$ and
$\Gamma(\tau_n^+  (X))$ have the same rational homology.
\end{cor}

\begin{proof}
Rationally any suspension is a wedge of spheres, so that $\Si^n X$
is rationally a $(n+1)$-fold suspension $\Si^{n+1} Y$, but
$\Lambda^n \Sigma^{n+1} Y \simeq \Gamma(\tau_n^+ ( \Si Y))$ has
the same homology as $\Gamma(\tau_n^+ ( X))$ by Corollary
\ref{corsec}.
\end{proof}

The corollary is not true in general, as the following
example shows.

\begin{ex} \label{not0}
The spaces
 $\Lambda^2(S^2 \vee S^2)$ and
$\Gamma(\tau_2^+ \sm (S^0 \vee S^0)  )$ have
not the same rational homology componentwise.

\end{ex}

\begin{proof}
Both spaces have components indexed by $\Z \x \Z$. The
bifiltration of the section space $\Gamma=\Gamma(\tau_2^+ \sm (S^0
\vee S^0) )$ by fiberwise bidegree is compatible with the
bifiltration, by number of particles, of the configuration space
$C(S^2,S^0 \vee S^0)$ of bicoloured particles on the sphere.
Compare also Example \ref{massa}. The bifiltration of the mapping
space $M=\Lambda^2(S^2 \vee S^2)$  is given by ordinary bidegree.

Let us denote the generators of $\tilde{H}_0(S^0 \vee S^0)
\subset \tilde{H}_0(\Om^2(S^2 \vee S^2)) $
by $x$ and $y$.
It turns out from Corollary \ref{corsec} that
the homology group $H_1(\Gamma_{m,n})$ is the quotient
of
$$H_1(\Om^2_{m,n}(S^2 \vee S^2))= \Q \{ [x,x]x^{m-2}y^n , [x,y]x^{m-1}y^{n-1}, [y,y]x^m y^{n-2}   \}$$
by
the subspace generated by $[x,x^{m-1}y^n ]$ and $[y,x^m y^{n-1}]$.
But these elements have respectively coordinates
$(m-1,n,0)$ and $(0,m,n-1)$
with respect to the basis above.
It follows that $H_1(\Gamma_{i,j}) \cong \Q^2$ for
$(i,j)=(1,1)$ or $(i,j)=(0,1)$ or $(i,j)=(1,0)$ and
$H_1(\Gamma_{i,j}) \cong \Q$ in all other cases.

On the other hand by Corollary \ref{cormaps} $H_1(M_{m,n})$ is the quotient of
$H_1(\Om^2_{m,n}(S^n \vee S^n))$ by the subspace generated by
$[x,x^{m-1}y^n ]+[x,x]x^{m-2}y^n$ and
$[y,x^m y^{n-1}]+[y,y]x^m y^{n-2}$, that have coordinates
$(m,n,0)$ and $(0,m,n)$. Thus
$H_1(M_{0,0})\cong \Q^3$ and $H_1(M_{i,j})\cong \Q$ for
any $(i,j) \neq (0,0)$.
In particular $M_{0,0}$ differs in homology
by all components of $\Gamma$.

\end{proof}

The next example shows that for $n$  even and at odd primes
the homology of the mapping space
$\La^n \Si^n X$
 does not
depend just on the additive homology of $X$, so that in particular
the section space and the mapping space have not the same homology.

\begin{ex} \label{mod3}

The spaces $X=S^2 \vee S^4 \vee S^6$ and $Y=\cp^3$ have the same
additive homology , but
$H_9(\La^2 \Si^2 X;\Z_3) \cong \Z_3 \oplus \Z_3$
and
$H_9(\La^2 \Si^2 Y;\Z_3) \cong \Z_3$.

\end{ex}

\begin{proof}
Let use denote the homology generators of the homology of
$X=S^2 \vee S^4 \vee S^6$ by $a_1,a_2,a_3$.
We must compute the Browder action mod 3
$$\lambda_i:\bigoplus_j \Si \tilde{H}_j(X) \ot H_{i-1-j}(\Om^2 \Si^2 X) \to H_i(\Om^2 \Si^2 X)$$
 for
$i=8,9$.  For $i=8$ the left hand side is generated by
$a_1 \ot [a_1,a_1]$ so that $\lambda_8 =0$ by
 Theorem 1.2 (6) in \cite{Cohen}.
For $i=9$ the left hand side has dimension 6, with generators
$a_1 \ot a_3, a_1 \ot a_1^3, a_1 \ot a_1 a_2,
a_2 \ot a_2, a_2 \ot a_1^2, a_3 \ot a_1$, and the right hand
side has dimension 5, with generators
$[a_1,a_3], [a_1,a_1]a_2, [a_1,a_2]a_1, [a_2,a_2], [a_1,a_1]a_1^2$.

By the Poisson relation (Theorem 1.2 (5) in \cite{Cohen})
$[a_1,a_1^3]=3[a_1,a_1]a_1^2 =0$ mod 3. Thus the image of
$\lambda_9$ has dimension 4, with generators
$[a_1,a_3],[a_1,a_1]a_2 + [a_1,a_2]a_1, [a_2,a_2], 2[a_2,a_1]a_1$,
and $H_9(\La^2 \Si^2 (S^2 \vee S^4 \vee S^6);\Z_3) \cong
ker(\lambda_8) \oplus coker(\lambda_9) \cong \Z_3 \oplus \Z_3$.

Let us consider now the case of the projective space $\cp^3$, with
homology generators $e_1, e_2, e_3$.
We must check the twisted Browder action
$$k_i:\bigoplus_j  \Si \tilde{H}_j(\cp^3) \ot
H_{i-1-j}( \Om^2 \Si^2 \cp^3) \to H_i(\Om^2 \Si^2 \cp^3)$$
for $i=8,9$.
By the proof of Corollary \ref{cormaps}
$k_*(\Si e_1 \ot y)=[e_1,y], \;
k_*(\Si e_2 \ot y)=[e_2,y] + [e_1,e_1]y$ and
$$k_*(\Si e_3 \ot y) = [e_3,y]+2[e_1,e_2]y-[e_1,e_1]e_1y.  $$
For the same reason as above $k_8=0$.
Moreover the image of $k_9$ is generated by the five elements
$[e_1,e_3],[e_1,e_1]e_2 + [e_1,e_2]e_1,
[e_1,e_1]e_2 + [e_2,e_2],
2[e_1,e_2]e_1 + [e_1,e_1]e_1^2 ,
[e_1,e_3]+2[e_1,e_2]e_1 - [e_1,e_1]e_1^2$.

Elementary linear algebra shows that they are linearly independent.
Thus $k_9$ is surjective, and
$H_9(\La^2 \Si^2 \cp^3;\Z_3)
\cong ker(k_8) \oplus coker(k_9) \cong \Z_3 $.

\end{proof}

We compute next the homology
of $\La^n S^{n+k}$ for $k \geq 0$.
We recall the computation of the homology of $\Om^n S^{n+k}$
for $n >1$.

Let $X$ be a space acted on by the operad of
 little $n$-discs. For example $X=\Omega^n Y$ or
$X=C(\R^n,Y)$.
We use the lower index notation for Dyer-Lashof operations
and do not consider $Q_0$ here.

\

a) Case $p=2$:

\smallskip

There are operations $Q_i:H_q(X;\Z_2) \to H_{2q+i}(X;\Z_2)$
for $0<i<n$.
(These operation are also named $\xi_i$ in \cite{Cohen}).
A string $Q_I$ of symbols $Q_i$ represents the composite
operation.
It is {\em admissible} if the sequence of indices is weakly
monotone. The empty string is admissible and represents the
identity.

\

b) Case $p$ odd:

\smallskip

There are operations
$Q_i(x):H_q(X;\Z_p) \to H_{pq+i(p-1)}(X;\Z_p)$, defined for
 $0<i<n$ and when
$i$ and $q$ have the same parity.
Note that
the operations $Q_i$ keep the parity and
the Bockstein operator $\beta$ switches parity.
An ordered string $Q_I$ of
symbols $Q_i$ and $\beta$
represents the composite operation.
We say that the operation $Q_I(x)$ is {\em admissible} if

1) the last symbol of the
non-empty string $Q_I$ is $Q_i$, and
$i$ has the same parity as the degree of $x$;

2) the sequence of indices is weakly monotone;

3) the indices of two adjacent $Q's$ have the same parity;

4) the indices of two
$Q's$ separated by a $\beta$ have opposite parity;

5) two adjacent $\beta$ do not appear.

\smallskip

In particular the empty string is admissible and represents
the identity.

\

c) In characteristic 0 only the identity is admissible.

\

Recall that $C(\R^n,S^k) \simeq \Om^n S^{n+k}$ for $k>0$ \cite{Bodig}.
Let $\io \in H_k(C(\R^n,S^k))$ be the fundamental class
of $S^k \simeq C_1(\R^n,S^k)$.

\begin{prop} {\rm \cite{Cohen}} \label{dyla}

1) The homology of $C(\R^n,S^k)$ is the
free commutative graded algebra on the admissible operations
$Q_I(\io)$,
and in addition on the admissible operations
$Q_J([\io,\io])$ if $n$ and $k$ have the same parity and
$char(K) \neq 2$.

2) The homology of $\Omega^n S^n$ is obtained
from $C(\R^n)=C(\R^n,S^0)$ by inverting
formally $\io$.

3) The basis of $H(C(\R^n))$ consisting of products
of admissible operations
 has an integral degree $\nu$ that keeps
track of the component, corresponding
to the number of particles. This degree is generated by
the rules $\nu(1)=0, \, \nu(\io)=1,\, \nu([\io,\io])=2,\, \nu(xy)=\nu(x)+\nu(y)$ and
$\nu(Q_i(x))=\nu(\beta Q_i(x))= p\nu(x)$.

\end{prop}

Now we are ready to state the computation.
The next theorem for $k>0$ gives also the homology of
$\La^n S^{n+k} \simeq C(S^n,S^k)$ by Proposition \ref{equal}(1).
Let $F$ denote the graded commutative algebra
generated by a set of vectors, and let $Q$ be
the set of all
admissible operations on $\io$ and $[\io,\io]$ except the identity.

\begin{thm} \label{basisconf}

\

1) For $n+k$ odd or for $char(K)=2$

$$H(C(S^n,S^k)) \cong H(C(\R^n,S^k)) \oplus
\Si^{n+k}H(C(\R^n,S^k)). $$

2) For $n$ odd, $k$ odd and $char(K) \neq 2$
$$H(\La^n S^{n+k}) \cong
H(\Om^n S^{n+k})/([\iota,\iota])
\oplus \Si^{n+k} H(\Om^n S^{n+k})/(\iota).$$

3) For $n$ even , $k$ even and
$p=char(K)$ odd

$$H(C(S^n,S^k)) \cong F(\io,Q) \oplus \io^{p-1} [\io,\io]
F(\io^p,Q) \oplus \Si^{n+k} [\io,\io] F(\io,Q) \oplus \Si^{n+k}
F(\io^p,Q).$$

4) For $n$ even, $k$ even and $char(K)=0$
$$H(C(S^n,S^k)) \cong F(\io) \oplus \Si^{n+k}[\io,\io]F(\io)
\oplus \Si^{n+k}\{1\} .$$

5) In the case $k=0$ the component degree of
$C(S^n)$ is obtained from the degree of
$C(\R^n)$ by the additional rule
$\nu(\Si^{n+k}x)=\nu(x)+1$.

\end{thm}

\begin{proof}

We apply Corollary \ref{corsec}.
The Browder action is given by $[\io,\_]$.
In case 1) the bracket $[\io,\io]$ is trivial.
In all cases
the bracket $[\io,x]$ is trivial if $x$ is any generator
of the algebra $H(\Om^n S^{n+k})$
other than $\io$.
Namely  $[\io,[\io,\io]]=0$ by 1.2(6) of \cite{Cohen},
and consequently
$[\io,Q_I(\io)]=0$ and $[\io,Q_J([\io,\io])]=0$ by 1.2(8)
and 1.3(9) of
\cite{Cohen}.
The action is a derivation with respect to the product
by 1.2(5) of \cite{Cohen}.
This completes the computation.

\end{proof}

The computation of the homology of configuration spaces
of even
spheres at odd primes is not deducible from
the methods of \cite{BCT}.

It remains to compute the homology of $\La^n S^n$.

\begin{thm} \label{basismaps}

\

\noindent 1) For $n$ odd or $char(K)=2 \quad H(\La^n S^n)=H(\Om^n
S^n) \ot H(S^n).$

\noindent 2) For $n$ even and $char(K)=0 \quad H(\La^n_0
S^n)=H(S^{n-1}) \ot  H(S^n)$ and $H(\La^n_k S^n) \cong
H(S^{2n-1})$ for  $k \neq 0$.

\noindent 3) For $n$ even and $char(K)=p$ odd $H(\La^n_k S^n)
\cong H(\Om^n_0 S^n) \ot H(S^n)$ if $p$ divides $k$ and

$H(\La^n_k S^n) \cong F(Q) \oplus \Si^{n} [\io,\io] F(Q)$
if $(p,k)=1$.

\end{thm}

\begin{proof}

Since the bundle $\tau_n^+$ is trivial, there is a homeomorphism
$\Gamma(\tau_n^+) \cong \La^n S^n$, and we can use either
Corollary \ref{corsec} or Corollary \ref{cormaps}. For $n$ odd
both the untwisted and the twisted Browder actions are trivial.
For $n$ even note that the base point of $\Gamma(\tau_n^+)$, the
section at infinity, is the antipodal map of $S^n$, that has
degree $-1$. If $\Gamma_k$ is the component of $\Gamma(\tau_n^+)$
of fiberwise degree $k$, then $\Gamma_{k+1} \cong \La^n_{k}S^n$.
This is compatible with corollaries \ref{corsec} and
\ref{cormaps}: the Browder action  $x \mapsto [\io,x]$ from
$\Om^n_k S^n$ to $\Om_{k+1}^n S^n$ computes the homology of
$\Gamma_{k+1}$, but the twisted Browder action from
$\Om^n_{k-1}S^n$ to $\Om^n_k S^n$
$\tau(x)=[\io,x]+[\io,\io]\io^{-1}x$ computes the homology of
$\La^n_k S^n$. They are compatible because $\io
\tau(x)=[\io,x]\io+[\io,\io]x= [\io,\io x]$. Cases 1) and 2) are
easy. In case 3) the multiplication by $\io^p$ induces an
isomorphism $H(\La^n_k S^n) \cong H(\La^n_{k+p}S^n)$. If $(p,k)=1$
then the composition by a degree $k$ map induces an isomorphism
$H(\La^n_1 S^n) \cong H(\La^n_k S^n)$. One concludes the
computation by keeping track of the component degree.

\end{proof}

\section{Homotopy pullbacks} \label{6}

Another approach to our spaces is given by homotopy pullbacks,
by decomposing the sphere $S^n$,
basis of the evaluation fibration, as union of two $n$-discs.

\begin{prop} \label{homaps}

The mapping space is the homotopy pullback
$$ \begin{diagram}
\Lambda^{n}\Si^nX & \rTo & \Si^nX \\
\dTo         &       & \dTo^c \\
\Si^nX           & \rTo^c   &\Lambda^{n-1} \Si^nX ,
\end{diagram}$$
where $c$
is adjoint to the projection
$S^{n-1}_+ \sm \Si^n X \to \Si^n X$.

\end{prop}

\begin{prop} \label{hosec}
The section space $\Gamma(\tau_n^+ ( X))$ is the homotopy pullback

$$ \begin{diagram}
\Gamma(\tau_n^+  (X)) & \rTo & \Si^nX \\
\dTo         &       & \dTo^c \\
\Si^nX           & \rTo^v   &\Lambda^{n-1} \Si^nX ,
\end{diagram}$$

where $v:\Si^n X \to \La^{n-1}\Si^n X$ is
adjoint to
$w:S^{n-1}_+ \sm X \to \Om^n \Si^n X$.
\end{prop}

\

The section and mapping spaces are both
total spaces of fibrations
 with same fiber and basis.
 The mapping space evaluation fibration admits
the section given by
 constant loops.
 Therefore for $i>0$
 $\pi_i(\La^n \Si^n X)=\pi_{i+n}(\Si^n X) \oplus \pi_i(\Si^n X)$.
The section space fibration admits just a stable section, given by
the inclusion $S^n_+ \sm X =C_1(S^n,X) \to \Gamma(\tau_n^+ ( X))$,
since $\Si(S^n_+ \sm X) \simeq \Si(\Si^n X \vee X)$.

\begin{ex} \label{last}
$C(S^2,\cp^2)$ and
 $\La^2 \Si^2 \cp^2$
 are not rationally homotopy equivalent,
 but they have the same homology over any field.
\end{ex}

\begin{proof}

We show first that the $\pi_5's$ of the two spaces have different ranks.
Rationally $\Si^2 \cp^2 \simeq S^4 \vee S^6$, so that
$rk(\pi_5(\La^2 \Si^2 \cp^2)) = rk(\pi_7(S^4 \vee S^6))
+ rk(\pi_5(S^4 \vee S^6))=1$.

The fiber square of Proposition \ref{hosec}
 induces a long sequence of rational
homotopy groups
$$\dots \pi_i(C(S^2,\cp^2)) \rTo \pi_i(\Si^2\cp^2)\oplus \pi_i(\Si^2 \cp^2)
\rTo^{c_\# \oplus v_\#}
\pi_i(\La \Si^2  \cp^2) \rTo \pi_{i-1}(C(S^2,\cp^2)) \dots$$

Now $\pi_5(S^4 \vee S^6)$ has rank 0, so that
rationally
$\pi_5(C(S^2,\cp^2))$ is the quotient of
$\pi_6(\La \Si^2 \cp^2)=\Q \{ s, u      \}$
by the image of $c_\# \oplus v_\#$, defined on
$\pi_6(\Si^2 \cp^2)^2 = \Q\{ s_1,s_2  \} $.

Here $s_1$ and $s_2$ represent the generators
coming from the wedge summand $S^6 \to S^4 \vee S^6$
and
$s=c_\#(s_1)=c_\#(s_2)$. Moreover $u$ is the composition
$S^6 \to \Omega S^4 \to \La S^4 \to  \La (S^4 \vee S^6)$,
with the first map adjoint        to the
Whitehead product $[\io_4,\io_4]:S^7 \to S^4$.

The group  $[S^1_+ \sm \cp^2 , \Om^2 \Si^2 \cp^2]$ splits
as $[\Si \cp^2,\Om^2 \Si^2 \cp^2 ] \oplus [\cp^2, \Om^2 \Si^2 \cp^2]$,
because after two suspensions $S^1_+$ and $S^1 \vee S^0$ are
equivalent as co-H spaces.
Rationally the first summand $\pi_5(\Om^2 \Si^2 \cp^2)=\Q\{\bar{u}\}$
has rank 1 and is detected
by the action on $\pi_5 \cong H_5$.

But $w: S^1_+ \sm \cp^2 \to \Om^2 \Si^2 \cp^2$ induces
the nontrivial action
$w_*([S^1] \ot e_2) = [e_1,e_1]$
where $e_i \in H_{2i}(\cp^2)$ is the generator,
by the proof of Corollary \ref{cormaps}.
Thus  $w$ corresponds in the decomposition above
to $\bar{u} \oplus i$, with $i$ adjoint to the identity.
By taking adjoints,
$v: \Si^2 \cp^2 \to \Lambda \Si^2 \cp^2$ represents the sum
in the group $[\Si^2 \cp^2, \Lambda \Si^2 \cp^2]$
of the injection $c: \Si^2 \cp^2 \to \La \Si^2 \cp^2$
and the composition
$\Si^2 \cp^2 \to S^6 \rTo^u \La \Si^2 \cp^2$.

It follows that $v_\#(s_2)=s + u$, and
$\pi_5(C(S^2,\cp^2))$ has rank 0.

The second assertion follows from
Corollary \ref{mod2} and the fact that away from the prime 2
$\Si^2 \cp^2 \simeq S^4 \vee S^6$.

\end{proof}

This shows that the rational homotopy type of
$C(S^2,X)$ does not depend just on the rational homology of
$X$, by taking $X=\cp^2$ and $X=S^2 \vee S^4$.

\

\noindent
Dipartimento di Matematica,     \\
Universit\`a di Roma ``Tor Vergata'', \\
Via della Ricerca Scientifica 1,  \\
00133 Roma, Italy \\
\\
e-mail: {\em salvator@mat.uniroma2.it}

\end{document}